\documentclass[12pt]{article}

\usepackage[margin=1in]{geometry}  
\usepackage{graphicx}              
\usepackage{amsmath}               
\usepackage{amsfonts}              
\usepackage{amsthm}                
\usepackage{multicol}
\usepackage{epsfig, epstopdf}
\usepackage{subfigure}
\usepackage{authblk}

\newtheorem{thm}{Theorem}[section]

\newcommand\T{\rule{0pt}{2.6ex}}

\begin{document}

\nocite{*}

\title{On parallel solution of ordinary differential equations}

\author[1]{Alejandra Gait\'an Montejo \thanks{ygaitanm@purdue.edu}}
\author[2]{Octavio A. Michel-Manzo \thanks{oalonso@math.cinvestav.mx}}
\author[2]{C\'esar A. Terrero-Escalante\thanks{Corresponding author: cterrero@ucol.mx}}
\affil[1]{Department of Mathematics, Purdue University, West Lafayette, IN 47907-2067, USA}
\affil[2]{Facultad de Ciencias, Universidad de Colima, Bernal D\'iaz 
del Castillo 340, Col. Villas San Sebasti\'an, Colima, Colima 28045, M\'exico}

\renewcommand\Authands{ and }

\maketitle

\begin{abstract}
  In this paper the performance of a parallel iterated Runge-Kutta method is  
compared versus those of the serial fouth order Runge-Kutta and Dormand-Prince 
methods. It was found that, typically, the runtime for the parallel method is 
comparable to  that of the serial versions, thought it uses considerably more 
computational resources. A new algorithm is proposed where full parallelization 
is used to estimate the  best stepsize for integration. It is shown 
that this new method outperforms the other, notably, in the integration of 
very large systems.
\end{abstract}

\section{Introduction} \label{sec:intro}

The numerical solution of an initial value problem given as a system of 
ordinary differential equations (ODEs) is often required in engineering and 
applied sciences, and is less common, but not unusual in pure sciences.
For precisely estimating asymptotic properties of the solutions, the global  
truncation errors must be kept lower than the desired tolerance during a very 
large number of iterations. This is usually achieved by using an adaptive 
algorithm for the estimation of the  largest step for integration yielding a 
local truncation error below the tolerance. Nevertheless, such a correction 
usually leads to a drastic increase of the computational time. On the other 
hand, the use of spectral methods to solve parabolic and hyperbolic partial 
differential equations (PDE) is becoming more and more popular. The spectral 
methods reduce these PDEs to a set of ODEs \cite{Boyd}. 
The higher the desired precision for 
the numerical solution, the larger the resulting system  of ODEs. 
Very large systems arise also in simulations of multi--agent systems \cite{Wooldridge}.
It can take hours to integrate this kind of systems over a few steps. 
Taking all the above into 
account, it can be concluded that devising improved algorithms to compute the 
numerical solution of ODE systems is still a very important task.

With the steady development of cheap multi-processors technology, it is 
reasonable to consider using parallel computing for speeding up real-time  
computations. Particularly interesting are the current options for small-scale 
parallelism with  a dozen or so relatively powerful processors. Several methods 
have been deviced with that aim (see for instance Refs. 
\cite{Burrage,Houwen,hairer}), many warranting a substantial reduction of the 
runtime. For instance, the authors in Ref.\cite{Houwen} claim that the 
performance of their  parallel method is comparable to that of the serial method 
developed by Dormand and Prince \cite{dop} and, in terms of the required number 
of evaluations of the right--hand side of the  ODEs, demonstrates a superior 
behaviour.

The aim of this work was twofold. Firstly, we wished to test these claims by 
solving ODEs systems with different degree of complexity over different ranges 
of time; secondly, we proposed and tested a new method which focuses on taking 
full advantage of parallel computing for estimating the optimal stepsize.
All our codes were written in C and for parallalel programing we used the 
OPENMP resources. The programs were tested in a server Supermicro A+ 1022GG-TF 
with 24 CPUs and 32 gigabytes of operational memory.

In the next section, we describe the numerical methods we used for our tests, 
the standard fourth order Runge-Kutta, a version of the Dormand-Prince method 
and the parallel iterated Runge-Kutta method proposed in Ref.\cite{Houwen}.
These last two methods are widely regarded to be amongst the best options for  
serial and parallel numerical solution of ODEs. It is also briefly described how 
the optimal stepsize is estimated in each case. In section \ref{sec:systems} the 
initial value problems used for testing the methods are described. Next, in 
section \ref{sec:tresults} we report the results of our comparison of the 
performance of these methods. In section \ref{sec:aspa} we introduce an adaptive 
stepsize parallel algorithm coupled to the Dormand-Prince integrator, and report 
the results of the corresponding tests. Finally, in section \ref{sec:concl} we 
present our conclusions.

\section{Numerical integrators and local error control} \label{sec:methods}

Let the initial value problem be specified as follows,
\begin{equation}
\dot{y} = f(t,y)\, , \quad y(t_0)=y_0\, .
\end{equation}
Here $y(t)$ is the vector solution at time $t$, dot stands for the derivative 
with respect to time and the right--hand side of the equation defines a 
vector field.

Our aim is to compare the performance of several methods for approximating  
the solution of this problem. All of them are members of the family of explicit 
Runge--Kutta methods and aproximate $y(t)$ at $t_{n+1}=t_{n}+h$ as
\begin{equation}
y_{n+1} = y_{n} + h\displaystyle\sum\limits_{i=1}^{s} b_{i}k_{i} \, ,
\end{equation}
where
\begin{eqnarray}
k_1 &=& f(t_n, y_n)\, , \nonumber \\
k_2 &=& f(t_n + c_2 h, y_n + h a_{21} k_1)\, , \nonumber \\
k_3 &=& f(t_n + c_3 h, y_n + h( a_{31} k_1 + a_{32} k_2)\, , \nonumber \\
\vdots \nonumber \\
k_s &=& f(t_n + c_s h, y_n + h( a_{s1} k_1 + a_{s2} k_2 + \cdots 
a_{s,s-1}k_{s-1})\, ,
\end{eqnarray}
and $s$ is known as the number of stages. Therefore, a method with $s$ stages 
usually requires, at least, $s$ evaluations of the right--hand side of the system 
at each iteration.

A Runge--Kutta method can be especified by a Butcher tableau like in table 
\ref{table:butcher}.
\begin{table}[h]
    \begin{center}
            $$\begin{array}{ c | c c c c c }
              0  \\
              c_{_2} & a_{_{21}}  \\
              c_{_3} & a_{_{31}} & a_{_{32}} \\
              \vdots & \vdots & \vdots & \ddots\\
              c_{_s} & a_{_{s1}} & a_{_{s2}} & \ldots & a_{_{s,s-1}}\\
              \hline
              \T& b_{_1} & b_{_2} & \ldots & b_{_{s-1}} & b_{_s}
            \end{array}$$
            \caption{Butcher tableau} \label{table:butcher}
    \end{center}
\end{table}
The order of a method is $p$ if the local truncation error is on the order of 
$\mathcal{O}(h^{p+1})$, while the total accumulated error is of order 
$\mathcal{O}(h^{p})$.
The higher the order of the method, the lower the error of the approximation, 
nevertheless, constructing higher order Runge--Kutta formulas is not an easy
task. To avoid increasing $s$ (and, therefore, the number of evaluations of $f$) 
a common alternative is to develope methods with adaptive stepsize.

For any numerical method, an estimate for the local truncation error while 
integrating from $t_n$ to $t_{n+1}=t_n+h$ is given by
\begin{equation}
\epsilon:=\|y_{n+1}-\bar{y}_{n+1}\|\, ,
\label{eq:error}
\end{equation}
where $\|\cdot\|$ stands for a given norm, and $y_{n+1}$ and $\bar{y}_{n+1}$ 
are the results of different numerical approximations of $y(t_{n+1})$. The 
stepsize yielding a local error below the tolerance ($Tol$) is then given by
\begin{equation}
h_{opt}=h\displaystyle\left(\frac{Tol}{\epsilon}\right)^{\frac{1}{p}}\, .
\end{equation}

\subsection{Fourth order Runge--Kutta method (RK4)}\label{sec:rk4}

The method given in table \ref{table:rk4} is the classical member of the family  
of Runge--Kutta methods.
\begin{table}[ht]
    \begin{center}
    \begin{tabular}{c|c c c c}
     0 & \\[.25cm]
    1/2 & 1/2 \\[.25cm]
    1/2 & 0 & 1/2 \\[.25cm]
    1 & 0 & 0 & 1 \\ \hline
    \T& 1/6 & 2/6 & 2/6 & 1/6
    \end{tabular}
    \caption{``The" Runge--Kutta method} \label{table:rk4}
    \end{center}
\end{table}

We used RK4 without a stepsize control mechanism.
Hence, in all our tests we choose the stepsize in such a way that the global  
error had the same order of those obtained by the methods with adaptive 
stepsize.

\subsection{Dormand-Prince method (DOP853)}\label{sec:dop}

In this method \cite{dop}, the last stage is evaluated at the same point as the  
first stage of the next step (this is the so-called FSAL property), so that the 
number of evaluations of $f$ is one less than the number of stages. Here there 
is no easy way to present the Butcher coefficients in a tableau, because it 
involves dealing with irrational quantities \cite{hairer}. The coefficients we 
use can be found in the code by  E. Hairer and G. Wanner  available in the site 
\cite{dopButcher}.

The approximations $y_{n+1}$ and $\bar{y}_{n+1}$ in equation (\ref{eq:error})  
correspond here to the results obtained using different orders, and the $k_i$ 
are determined by minimizing the error of the higher order result. As a matter 
of fact, in the version we use two comparisons are made, one between 8th and 5th 
orders, and the second one between 8th and 3th orders. Then, the error is 
estimated using \cite{hairer}:
\begin{equation}
\epsilon=\epsilon_5 \frac{\epsilon_5}{\sqrt{\epsilon_5^2+0.01\epsilon_3}} \, .
\end{equation}

\subsection{Parallel iterated Runge--Kutta method (PIRK10)}\label{sec:pirk}
    Let us consider a s-stage Runge--Kutta method given by the coefficients
    $$A=(a_{ij})_{_{i,j=1}}^{^s}, \qquad B^{^T}=(b_{_1},\ldots,b_s), \qquad  
      C=(c_{_1},\ldots,c_s)^{^T}$$
    and let $y_1$ be defined as:
        \begin{flalign}
        k_{i}^{^{(0)}}&=f(x_{_0},y_{_0}) \nonumber\\
        k_{i}^{^{(\ell)}}&=f(x_{_0}+c_{i}h,y_{_0}+ 
		    h\displaystyle\sum\limits_{j=1}^{s} 
		    a_{{ij}}k_{j}^{^{(\ell-1)}}) \ \ \ \ \  
		    \ell=1,\ldots,m\label{PIRK}\\
        y_1&=y_0+h\displaystyle\sum\limits_{i=1}^{s} b_{i}k_{i}^{^{(m)}} 
	    \nonumber
    \end{flalign}

    Here $m$ is the number of iterations used to estimate $k_{i}$.
As it is shown in \cite{Houwen}, provided that $s$ processors are  available, 
this scheme represents an explicit Runge--Kutta method. Furthermore, since each 
$k_{i}^{^{(\ell_{_0})}}$ can be computed in parallel, 
we have the following theorem,
    \begin{thm}
     Let $\{A,B^T,C\}$ define an s-stage Runge--Kutta method of order $p_0$. Then 
the method defined by \eqref{PIRK} represents an $(m+1)-stage$ explicit  
Runge--Kutta method of order $p$, where $$p=\min\{p_0,m+1\}.$$
\label{thm:thm}
    \end{thm}
    One of the advantages of this method is that if we set $m=p_0-1$, then the  
order of the method is equal to the number of stages, which results in less 
right--hand side evaluations (sequentially). In general the number of stages of 
explicit Runge--Kutta methods is greater than the order of the method, therefore 
if an explicit method is used the required number of processors is greater as 
well. 

Along the lines in Ref.\cite{Houwen} and with the Butcher coefficients in table 
\ref{table:RK10} in the Appendix, we implemented a parallel iterated Runge--Kutta 
method of order $10$. Here, $y_n$ and $\bar{y}_n$ in equation (\ref{eq:error}) 
correspond to the results obtained using different number of iterations for 
approximating $k_{i}$, 
\begin{equation}
y_{_{n+1}}=y_{_n}+h\displaystyle\sum\limits_{i=1}^{s} b_{i}k_{i}^{^{(m)}}
\end{equation}
and
\begin{equation}
\bar{y}_{_{n+1}}=y_{_n}+h\displaystyle\sum\limits_{i=1}^{s} 
		  b_{i}k_{i}^{^{(m-1)}} \, .
\end{equation}

\section{Initial value problems} \label{sec:systems}

Next, we list and briefly describe the systems of ODEs that we have used to test 
the above numerical methods.

\subsection{Simple harmonic oscillator (HO)} \label{ssec:HO}

As a first initial value problem we chose:
    \begin{equation}
        \begin{cases}
        \dot{y}_{1}=y_2 &\qquad y_{1}(0)=0,\nonumber\\
        \dot{y}_{2}=-y_1 &\qquad y_2(0)=1\text{.} \label{eq:ho}
        \end{cases}
   \end{equation}

Since this system is readily integrable, we used it to assess the quality of the 
numerical results by comparing with the analytical ones.

\subsection{H\'enon-Heiles system (HH)} \label{ssec:HH}

This is a Hamiltonian system which describes the nonlinear dynamics of a star  
around a galactic center when the motion is constrained to a plane \cite{hh}:
    \begin{equation}
        \begin{cases}
        \dot{y}_{_1}&=y_{_{2}}\nonumber\\
        \dot{y}_{_2}&=-y_{_1}-2y_1y_3\\
        \dot{y}_{_3}&=y_{_4}\nonumber\\
        \dot{y}_{_4}&=-y_{_3}-y^2_{_1}+y^2_{_3}  \label{eq:HH}
        \end{cases}
   \end{equation}
Since the Hamiltonian $H$ is a constant of motion, it can be used to assess the  
precision of the numerical solution. We choose initial conditions such that 
$H=1/6$, yielding a chaotic solution.

\subsection{Replicated H\'enon-Heiles system (HH100)} \label{ssec:HH100}

To force the integrators to work a little bit harder we constructed a new system
by replicating the H\'enon-Heiles system $100$ times,
resulting in a nonlinear system with $400$ equations:
    \begin{equation}
        \begin{cases}
        \dot{y}_{_{4i+1}}&=y_{_{4i+2}}\nonumber\\
        \dot{y}_{_{4i+2}}&=-y_{_{4i+1}}-2y_{_{4i+1}}y_{_{4i+3}}\\
        \dot{y}_{_{4i+3}}&=y_{_{4i+4}}\nonumber\\
        \dot{y}_{_{4i+4}}&=-y_{_{4i+3}}-y^2_{_{4i+1}}+y^2_{_{4i+3}}  \, ,
        \end{cases}
\label{eq:HH100}
   \end{equation}
with $i=0,1,\cdots,99$.

\subsection{Gravitational collapse in AdS (GC40) and (GC10)} \label{ssec:ads}

We also tested the methods by solving the system obtained from the Einstein 
field equations for the gravitational collapse of a scalar field in anti de 
Sitter spacetime \cite{deOliveira:2012ac}. Using the Galerkin method \cite{Boyd}
the 10 coupled hyperbolic-elliptic nonlinear partial differential equations were
converted to a set of $40$ nonlinear ordinary differential equations. The 
corresponding solutions were shown to be chaotic too \cite{deOliveira:2012ac}.

Finally, the last system we used was obtained by reducing the previous one to  
ten equations \footnote{Any of the equations of the these two last systems fills 
several pages. The systems in C code are available from the authors.}.

\section{Tests results} \label{sec:tresults}

To test the methods we ask for the numerical solution of the corresponding problem  
starting from $t_0$ and up to a given $t_{end}$, such that the straigthforward 
integration with step $h_0=t_{end}-t_0$ yields a result with an error above the 
desired tolerance. This implies that, typically, a number of intermediate 
integrations will be required. 

In table \ref{table:smalltime} is shown the order of the runtime in seconds 
taken for solving the HO and HH problems in the time interval $0\leq t \leq 
2000$ using RK4, DOP853 and PIRK10. In the methods with an adaptive stepsize 
algorithm we have used a tolerance of $10^{-15}$, what corresponded to using a 
step $h=0.01$ in the RK4.
  \begin{table}[ht]
   \begin{center}
      $$\begin{array}{c|c|c}
         \T & HO & HH\\[.2cm]\hline
	DOP853\T & 10^{-2} & 10^{-2}\\[.2cm]
	PIRK10 & 10^{-1}  & 10^{-1}\\[.2cm]
	RK4 & 10^{-1}  & 10^{-1}\\[.2cm]
	\hline
        \end{array}$$
	\caption{Order of the runtime for the HO and HH problems.} \label{table:smalltime}
   \end{center}
  \end{table}
In all the following tests the PRIK10 used its optimal number of 5 processors.

  As we can see very similar results were obtained with the three methods, and  
even though DOP853 seems to be faster, the differences are very small.
Nevertheless, the serial methods can be considered to be better than PIRK10 
because they are easier to implement and require significantly less 
computational resources for execution.

Searching for a bigger runtime difference we tested the HH100 problem keeping  
the same tolerance for DOP853 and PIRK10, but now in the time interval $0\leq t 
\leq 5000 $. This implied to use $h=0.001$ in the RK4. In this case the RK4 and 
PIRK10 recorded a runtime of $\sim206$ seconds and $\sim75$ seconds 
respectively, both greater than the $\sim 11$ seconds obtained with DOP853.

At this point we recall that, according with theorem \ref{thm:thm},
by using 5 processors the PRIK10 method at each timestep does 9 evaluations of  
the right-hand-side of the corresponding problem. This is to be contrasted with 
the, at least, 11 evaluations done at each timestep by the DOP853. Therefore, 
since according with the above results the serial method outperforms the 
parallel one, we conjecture that this due to a parallel overhead problem, i.e.,
the amount of time required to coordinate parallel tasks is larger than the time 
required for evaluating the system right--hand side.

To verify this conjecture we tested the methods with the huge system of problem  
GC40. In this case we integrated the system over the small time interval $0\leq 
t \leq 0.1$, with a tolerance of $10^{-6}$, what corresponded to using a step 
$h=0.0001$ in the RK4. The results are presented in table 
\ref{table:colapso1}.
    \begin{table}[h!]
    \begin{center}
        \begin{tabular}{c c}
        \hline
        \T& $Time$  \\[.2cm] \hline
         DOP853 \T& $> 6$ days \\[.25cm]
         PIRK10 & $\approx 6$ hrs \\[.25cm]
         RK4& $\approx 2$ days \\[.25cm]
         \hline
        \end{tabular}
        \caption{Gravitational collapse runtime.} \label{table:colapso1}
        \end{center}
    \end{table}
 We can observe that the performance of PIRK10 was way better than DOP853 and  
RK4, being DOP853 unable to solve the system after six days.

\section{Adaptive stepsize parallel algorithm (ASPA)}\label{sec:aspa}

Since parallelizing the integrator does not seems to be helpful,
we opted for a different approach, that is, to parallelize the choice of an optimal integration step.

Let us consider an embedded Runge--Kutta method, which allows us to estimate the local error $\epsilon$. 
Given an initial step $h_0$ and a tolerance $Tol$,
for integrating from $t_0$ to $t_{end}$ with $N_{CPU}$ processors, 
the next step is determined as follows:
  \begin{enumerate}
   \item Each processor $P_i$, with $i=1,\ldots,N_{CPU}$, integrates the
system from $t_n$ to $t_n+ih_n$ and estimates the local error $\epsilon_i$.
   \item $m=\max_i\{i\ |\ \epsilon_i\leq Tol\}\cup\{0\}$.
   \item
    $ 
h_{n+1}=\Big(\frac{2N_{CPU}-1}{N_{CPU}+1}m+\frac{N_{CPU}}{2N_{CPU}-1}\Big) 
\frac {h_n}{N_{CPU}+1}$.
  \item $t_{n+1}=t_n+mh_n$.
  \item All the above steps are repeated while $t_k<t_{end}$.
  \end{enumerate}

\begin{figure}[h]
\centering
\mbox{{\includegraphics[width=13cm, height=6cm]{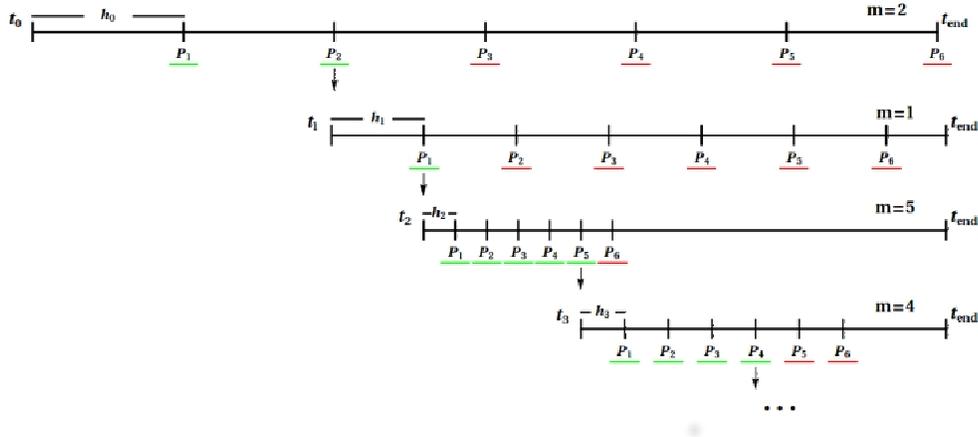}}}
\caption{An illustration of how the adaptive stepsize parallel algorithm 
could work with $N_{CPU}=6$ processors.}
\label{fig:aspa}
\end{figure}

Figure \ref{fig:aspa} is an illustration of how the stepsize could change with 
respect to $h_0$, depending on the number $m$ of processors yielding an 
acceptable result.

The interval of $N_{CPU}$ black vertical bars is the amount of time probed by 
the integration using the initial step $h_0$. 
In our computations $h_0$ is assumed to be the total length of integration over the number of processors. 
A green horizontal line below a processor label indicates a successful integration, 
otherwise, a red line is used.

\subsection{Designing the stepsize recurrence}

In a given iteration we define success as obtaining an integration result with a local error below the user defined tolerance.
The aim is to maximize the probability of success in each iteration, 
i.e., to determine the step $h_n$ such that it is obtained the biggest possible number of successful processors $m$
amongst the total number of available CPUs ($N_{CPU}$). 
We define $h_n$ as a function of $m$, 
keeping $N_{CPU}$ constant. 
So, if with a given integration step less than half of the processors are successful, 
then the next integration step needs to be smaller ($h_{n+1} < h_{n}$). 
Otherwise, we increase the integration step. 
This way, each integration becomes more efficient, both in amount of time and precision.
This idea is summarized with the following expresion:
\begin{align*}
 h_{n+1}&\approx \left(\frac{2}{N_{CPU}}m+\epsilon\right)h_n,\quad \text{for some } 
\epsilon>0,\\
 &=
 \begin{cases}
  \left[(1+\frac{2}{N_{CPU}} k) + \epsilon\right]h_n & \text{ if } m = N_{CPU}/2 + k,\\
  \left[(1-\frac{2}{N_{CPU}} k)+\epsilon\right]h_n & \text{ if } m = N_{CPU}/2 - k,\\
 \end{cases}\\
\end{align*}
for some integer $k\in[0,N_{CPU}/2]$.
We need $\epsilon$ to keep $h_{n+1}$ finite, 
even when $m=0$ and
it has to be less than one for the step to always decrease in this particular case.
A reasonable proposal is then, to carry the integration in half the interval when $m=0$,
i.e.,
\begin{equation}
  h_{n+1}\approx \frac{2}{N_{CPU}}mh_n+\frac{1}{2N_{CPU}}h_n\, . \nonumber
\end{equation}
On the other hand, 
if $m$ is large enough, the integration step will nearly doubles.
If, for instance, 
this happens sequentially, for typical initial value problems 
there is a high probability that the next $m$ will be very small,
making a poor use of the available CPUs.
To avoid this,
we finally propose the following recurrence:
\begin{equation}
  h_{n+1}=\Bigg[\frac{2N_{CPU}-1}{(N_{CPU}+1)^2}m
+ \frac{N_{CPU}}{(2N_{CPU}-1)(N_{CPU}+1)}\Bigg] h_n \, .
\label{eq:hn}
\end{equation}
Since $m$ is a function of $h_{n}$, this is a first order nonlinear map.
It yields, $h_{n+1} > h_n$ if
\[
 m>\dfrac{(N_{CPU}+1)(2N_{CPU}^2-1)}{(2N_{CPU}-1)^2} > \frac{N_{CPU}}2+1 \, ,
\]
and
$h_{n+1} < h_n$ if
\[
 m<\dfrac{(N_{CPU}+1)(2N_{CPU}^2-1)}{(2N_{CPU}-1)^2} < \frac{N_{CPU}}2 \, .
\]
Moreover, $m=0$ implies $1/2>h_{n+1}/h_{n}>0$, 
and when $m=N_{CPU}$, then $2>h_{n+1}/h_{n}>3/4$.
In consequence, this expression has the desired properties;
a large integration step will ultimately lead to a low $m$
that, in turn, will decrease the stepsize and, then, increase $m$. 
This way, we expect $m$ to converge to the optimal value $N_{CPU}/2$. 

Nevertheless, it is not desirable that the stepsize occurs to be insensitive to the given integration interval.
Thus, the map $h_{n+1}(h_n)$ was also designed to not have fixed points.
Note that requiring $m\leq N_{CPU}$, implies $N_{CPU}>2$,
i.e., there are not fixed points when using less than 3 CPUs.
For the remaining cases, $h_{n+1} = h_n$ give us the condition for the map to have fixed points:
\[
 m=\dfrac{(N_{CPU}+1)(2N_{CPU}^2-1)}{(2N_{CPU}-1)^2}\, .
\]
Let us prove that, whereas $m\in\mathbb{Z}$, the righ hand side of the above expresion is never an integer.
Suppose there is a $d\in\mathbb{Z}$ such that $d|(2N_{CPU}^2-1)$ and $d|(2N_{CPU}-1)$ 
\footnote{Here $d|f$ stands for $d$ divides $f$.}.
Since $2N_{CPU}^2-1 = (N_{CPU}+1)(2N_{CPU}-1)-N_{CPU}$ then, $d|N_{CPU}$.
So, by assumption $d|[(2N_{CPU}-1)-N_{CPU}]$, leading to $d|(N_{CPU}-1)$.
Therefore, considering that $d|N_{CPU}$, $d|(N_{CPU}-1)$ and $\gcd(N_{CPU},N_{CPU}-1)=1$,
we conclude that $d=1$.
In turn this implies
$\gcd(2N_{CPU}^2-1,2N_{CPU}-1)=1$,
and, this way, $(2N_{CPU}^2-1)/(2N_{CPU}-1)^2$ is not an integer.
Thus, $(N_{CPU}+1)(2N_{CPU}^2-1)/(2N_{CPU}-1)^2$ is an integer if and only if,
$(2N_{CPU}-1)^2|(N_{CPU}+1)$.
But, recalling that $N_{CPU}>2$, then $(N_{CPU}+1)/(2N_{CPU}-1)^2<1$
and we get a contradiction because by definition $m\in\mathbb Z$. 
Therefore $\{h_n\}$ has no fixed points.

\begin{figure}
\centering
\mbox{\subfigure{\includegraphics[width=3in, height=2in]{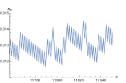}}\quad
\subfigure{\includegraphics[width=3in, height=2in]{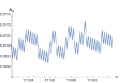}}
}
\mbox{\subfigure{\includegraphics[width=3in, height=2in]{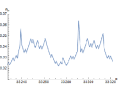}}\quad
\subfigure{\includegraphics[width=3in, height=2in]{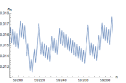}}
}
\mbox{\subfigure{\includegraphics[width=3in, height=2in]{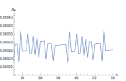}}\quad
\subfigure{\includegraphics[width=3in, height=2in]{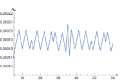}}
}
\caption{Oscillatory behaviour of $h_n$ vs $n$. 
Top: HH system for $N_{CPU}=10$ and $N_{CPU}=15$. Middle: HH100 system for 
$N_{CPU}=12$ and $N_{CPU}=14$. Right: Collapse (GC40) for 12 and 17 processors 
and $t=.1$.} 
\label{fig:hnbeh}
\end{figure}

Finally, note that while increasing the value of $N_{CPU}$, $h_0$ becomes smaller, 
implying a big number of integration steps in the beginning of the process. 
However at a given time, since there are not fixed points, $h_n$ should show a bounded oscillatory behaviour around the optimal stepsize.
It would imply that the proposed recurrence has an attractor,
i.e., asymptotically, the process of integration will settled down around an optimal stepsize independently of its initial value.
Indeed, this can be seen in figures \ref{fig:hnbeh} where we show some numerical realizations of $h_n(n)$
for different initial value problems and number of CPUs.

\subsection{Testing ASPA}

We tested the above described algorithm by coupling it to a version of the  
serial DOP853. We then compared the performances of the serial DOP853 and the 
DOP853 with ASPA (DOP853-ASPA). With this aim we calculated the difference of 
the number of stepsize corrections and the difference of runtime required to 
reach $t=t_{end}$ as function of the tolerance for a fixed number $N_{CPU}$ of 
processors. We also calculated the same differences but as function of the 
number of processors with the tolerance fixed to $10^{-15}$. The actual values 
of the runtime for each case are given in correponding tables in the appendix \ref{sec:runtimes}.
In figures \ref{fig:ASPAHH} the results for the HH problem are presented.
\begin{figure}
\centering
\mbox{\subfigure{\includegraphics[width=3in, height=3in]{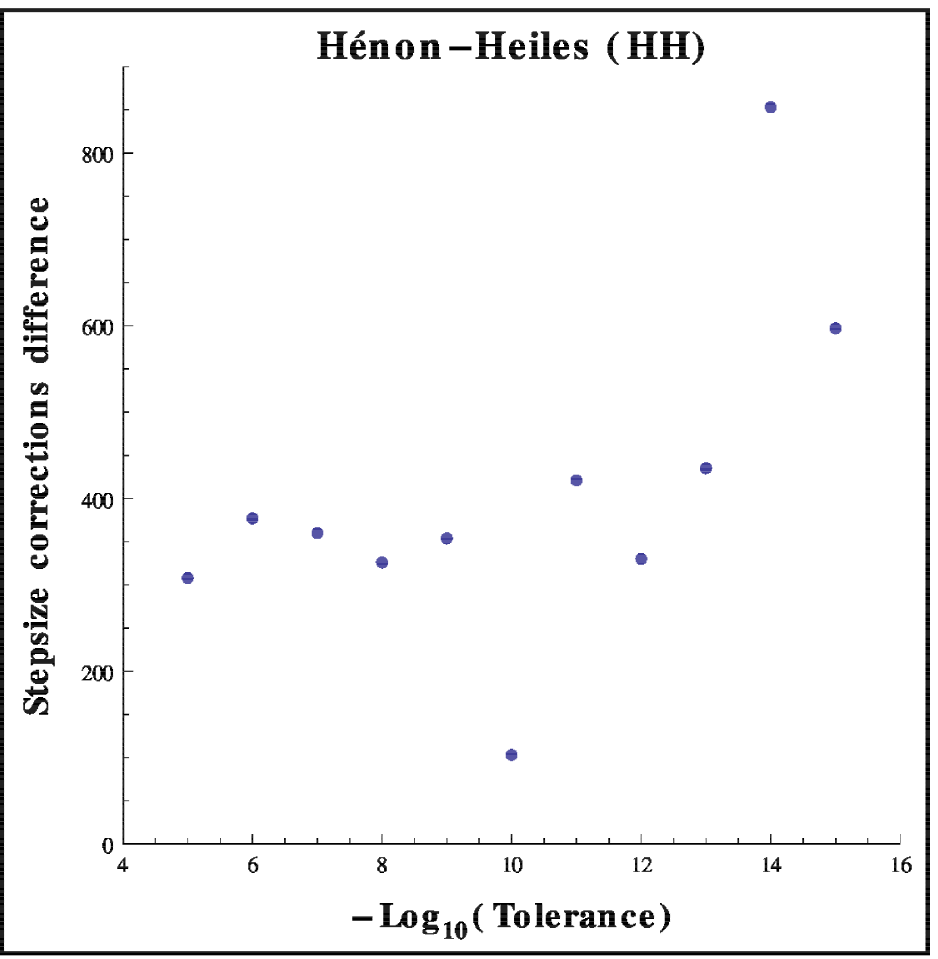}}\quad
\subfigure{\includegraphics[width=3in, height=3in]{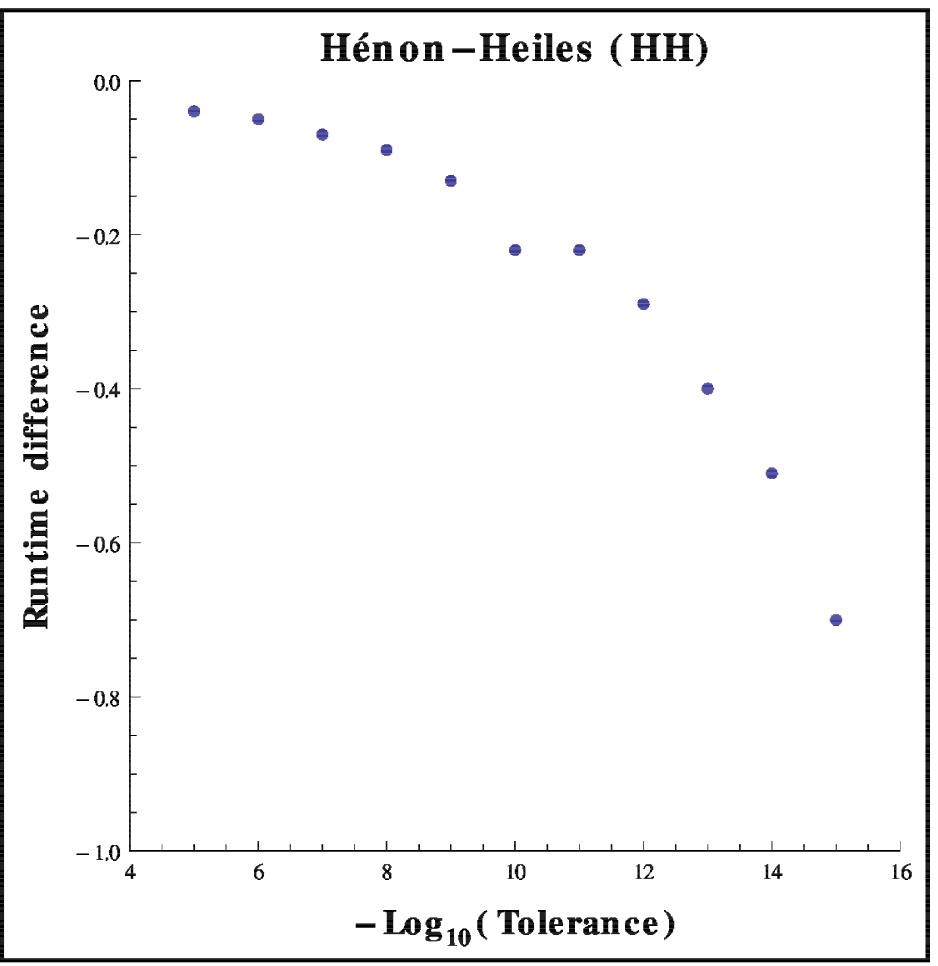} }}
\mbox{\subfigure{\includegraphics[width=3in, height=3in]{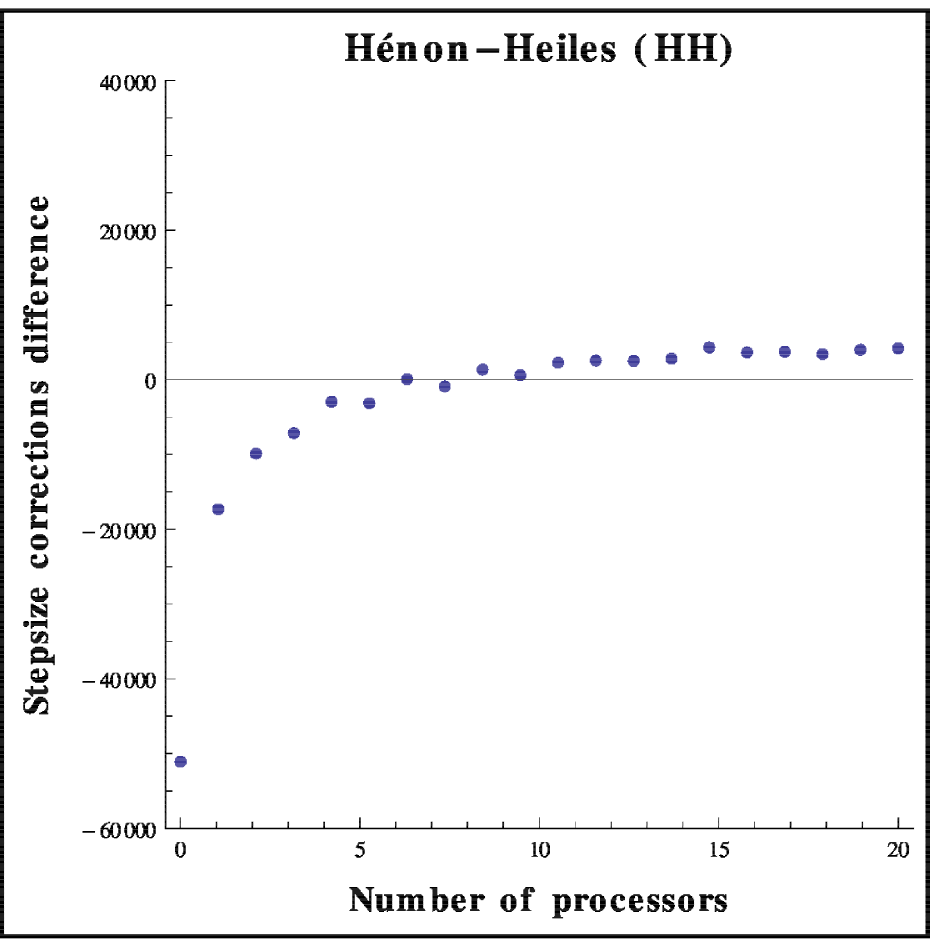}}\quad
\subfigure{\includegraphics[width=3in, height=3in]{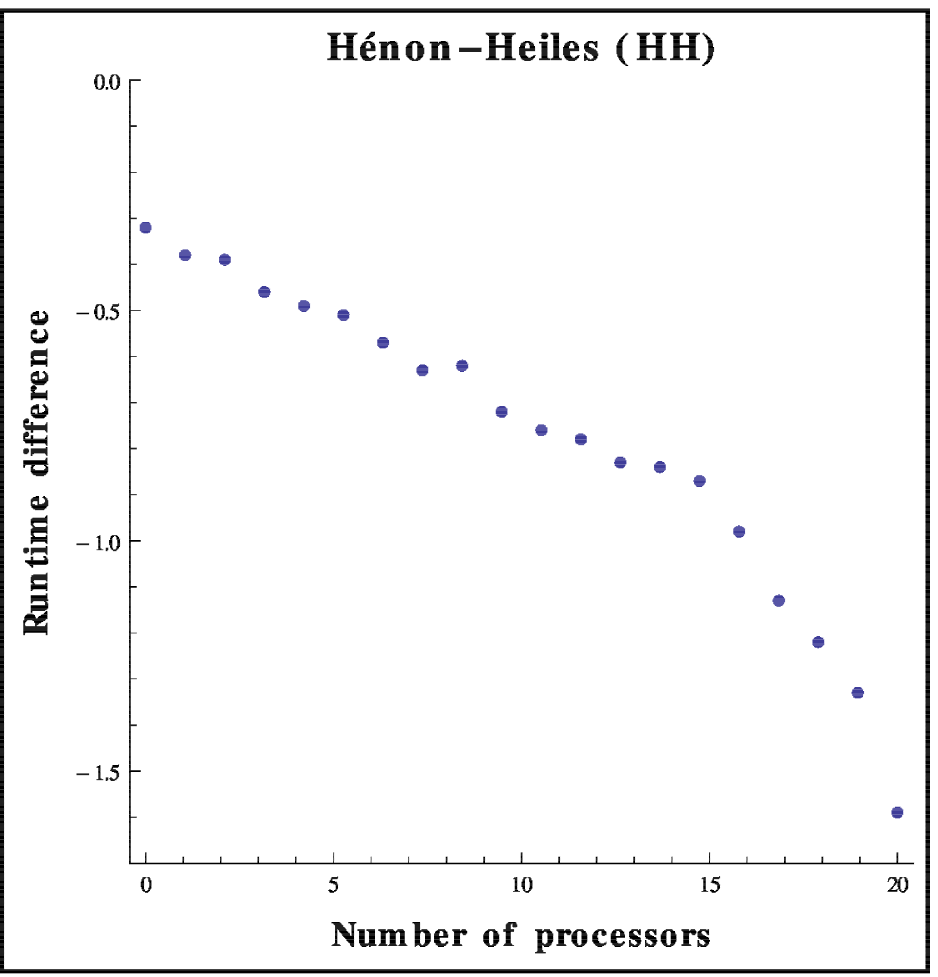} }}
\caption{DOP853 minus DOP853-ASPA for HH.
Left top: stepsize corrections vs tolerance. Left bottom: stepsize corrections  
vs number of processors. Right top: runtime vs tolerance. Right bottom: 
runtime vs number of processors.} \label{fig:ASPAHH}
\end{figure}
Here $t_{end}=5000$ and $N_{CPU}=10$.
In the two top panels we can see that, even if the DOP853 requires more 
stepsize corrections to reach the required tolerance, it does it in relatively 
less runtime. From the two bottom panels we draw the unexpected conclusion that
the runtimes are comparable only when the number of stepsize corrections  
required by DOP853-ASPA is significantly more than that required by DOP853.
This happens when using five or less processors. Moreover, notice that in the 
bottom panel the differences are all calculated with respect of the fixed 
number obtained with DOP853 (where  $N_{CPU}=1$). 
It means that, as expected, increasing $N_{CPU}$, the 
number of stepsize corrections in DOP853-ASPA decreases, nevertheless, the 
corresponding runtime increases. All these observations hint that, when more 
processors are used, at each iteration the parallel overhead is more important than 
the time required for integration.

\begin{figure}
\centering
\mbox{\subfigure{\includegraphics[width=3in, 
height=3in]{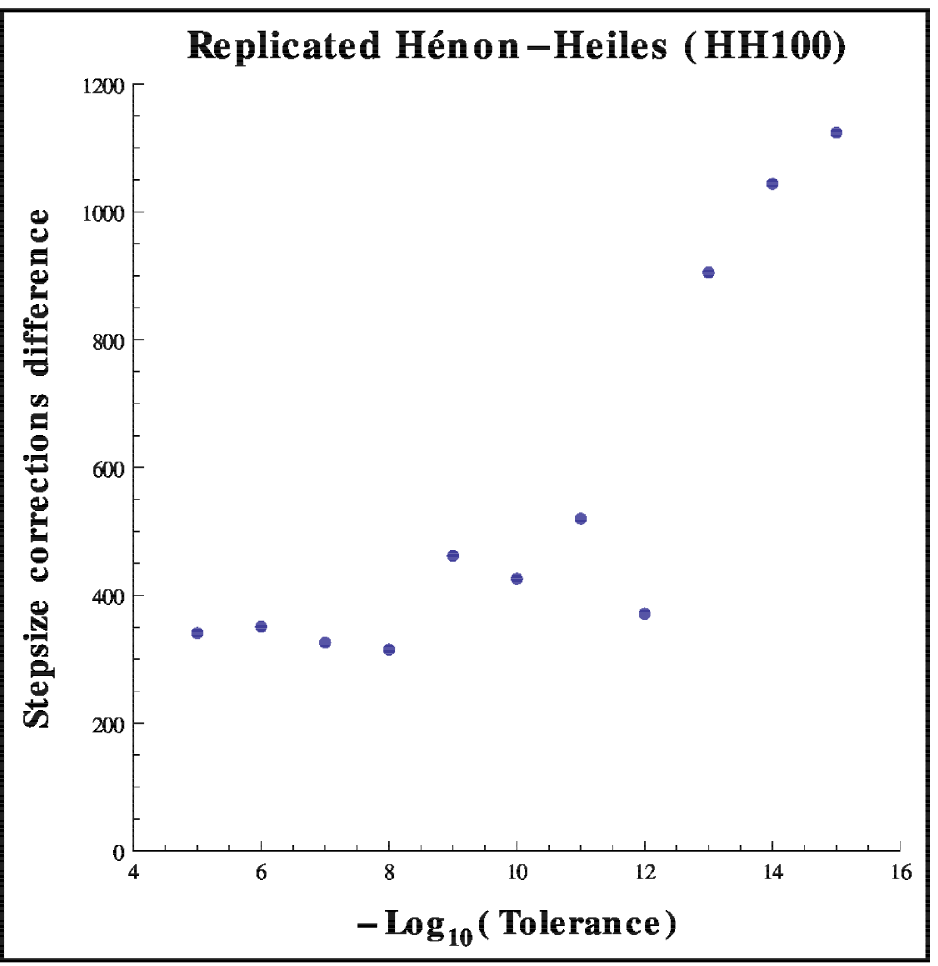}}\quad
\subfigure{\includegraphics[width=3in, height=3in]{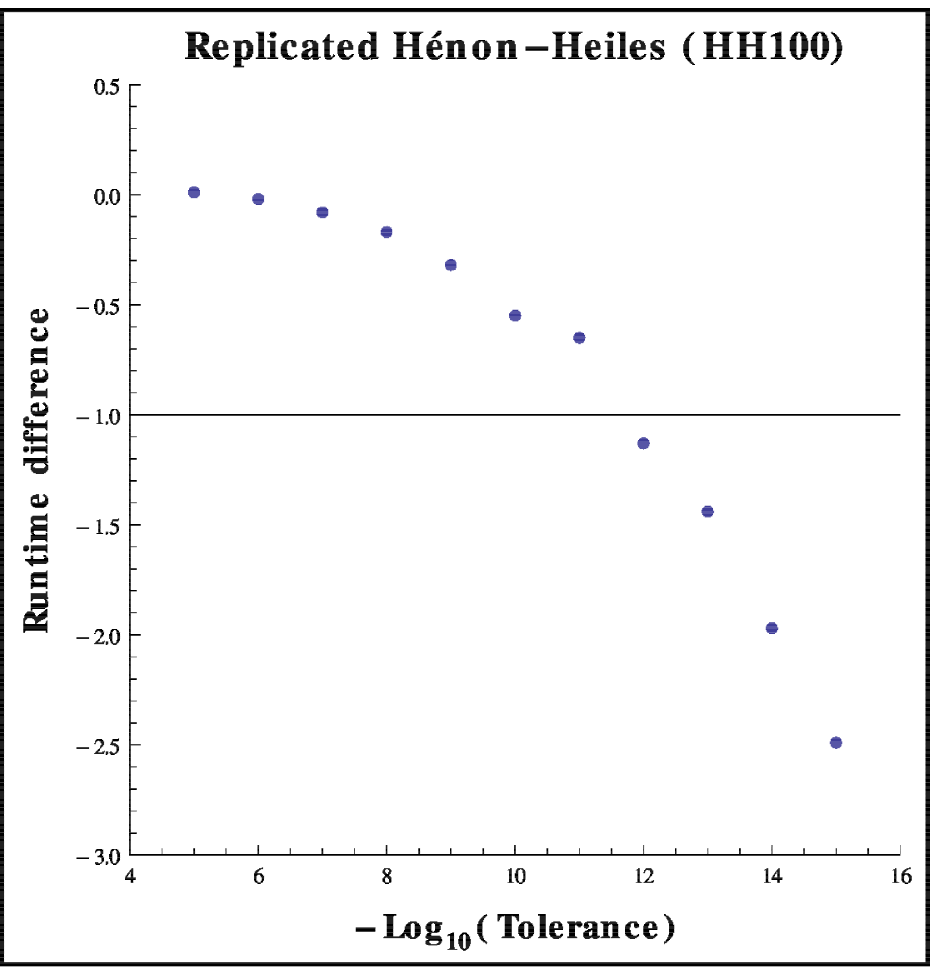}}
}
\mbox{\subfigure{\includegraphics[width=3in, 
height=3in]{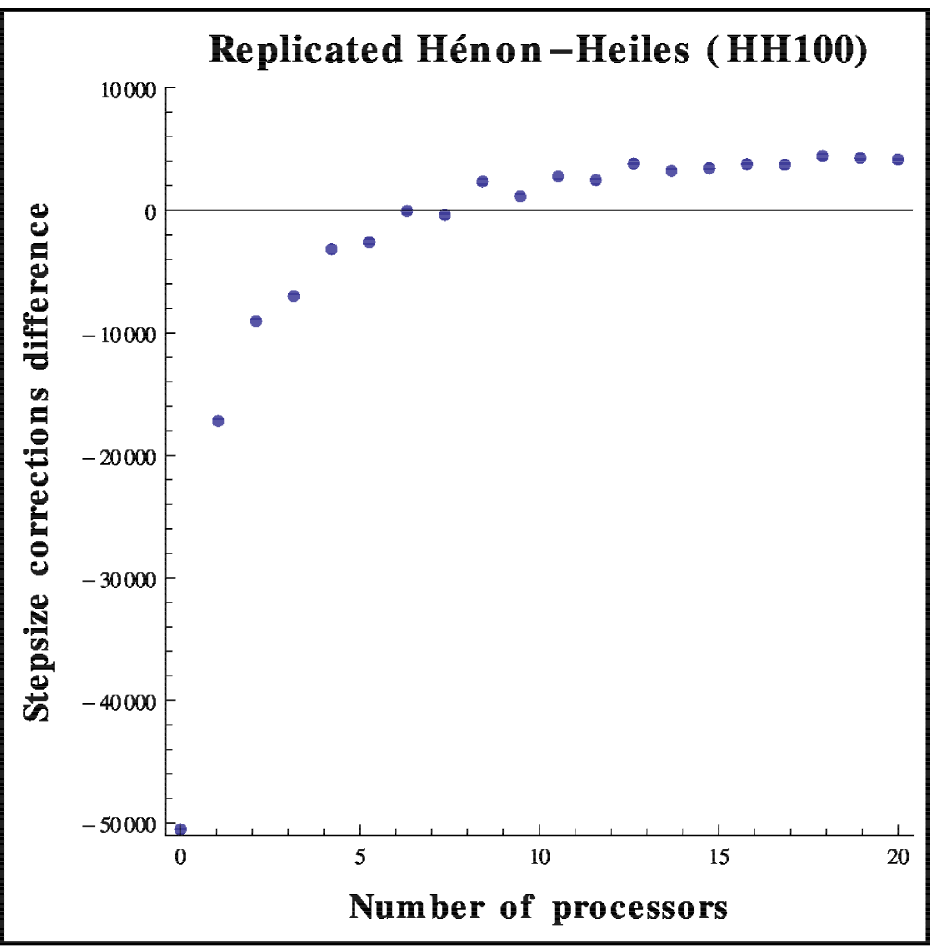}}\quad
\subfigure{\includegraphics[width=3in, height=3in]{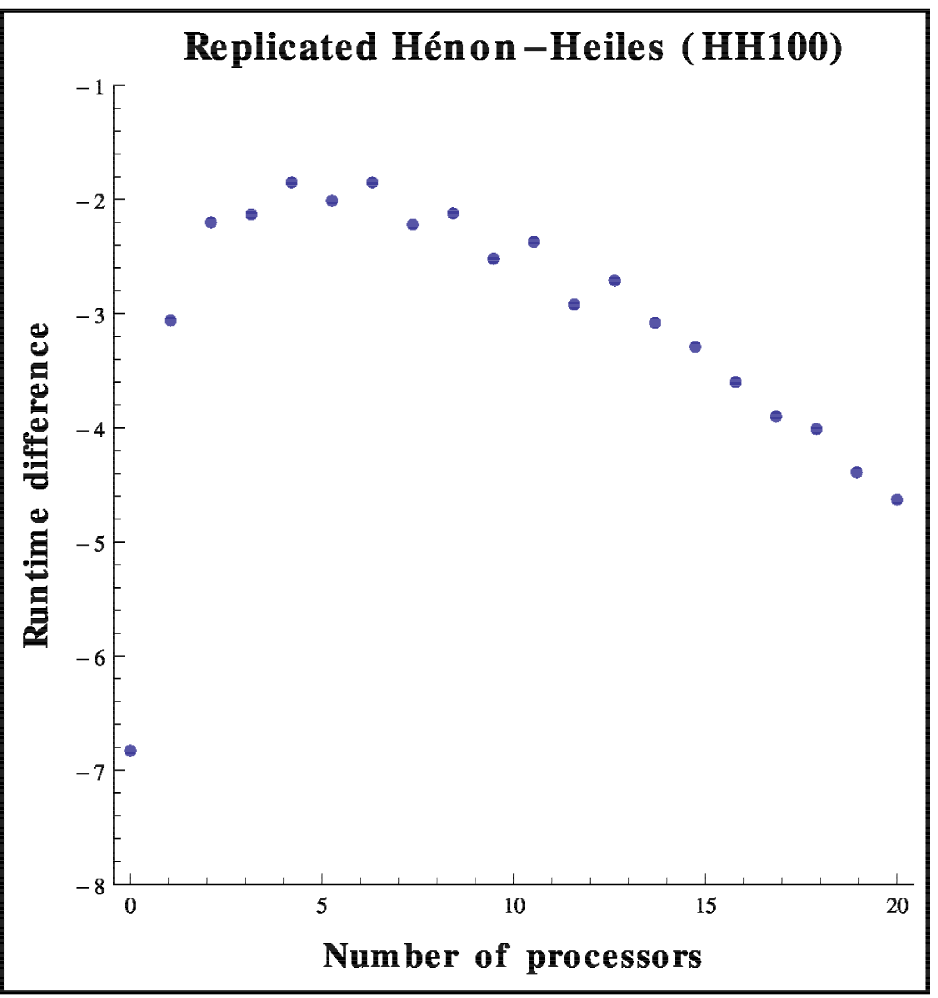} }
}
\caption{DOP853 minus DOP853-ASPA for HH100. Left top: stepsize corrections vs 
tolerance. Left bottom: stepsize corrections vs number of processors. Right 
top: 
runtime vs tolerance. Right bottom: runtime vs number of processors.} 
\label{fig:ASPAHH100}
\end{figure}

To determine whether this is the case, we tested the HH100 problem.
The results are presented in figures \ref{fig:ASPAHH100}.
Here $t_{end}=5000$ and $N_{CPU}=10$.
As in the case of HH, here DOP853 requires more stepsize corrections than  
DOP853-ASPA, nevertheless, for low tolerances the parallel algorithm performs 
slightly better than the serial. This could be due to the fact that for the 
HH100 problem the amount of time used for the evaluation of the RHS is 
comparable with the parallel overhead and that, for tolerances greater than 
$10^{-7}$, $N_{CPU}=10$ processors are good enough to probe the whole time 
interval up to $t_{end}=5000$ in very few stages.

Trying further to make the number of evaluations of the RHS to have a larger  
weight in the runtime, we tested the problem of the gravitational collapse, but 
the reduced version GC10, because the DOP853 was able to integrate it in a 
reasonable runtime.

\begin{figure}
\centering
\mbox{\subfigure{\includegraphics[width=3in, height=3in]{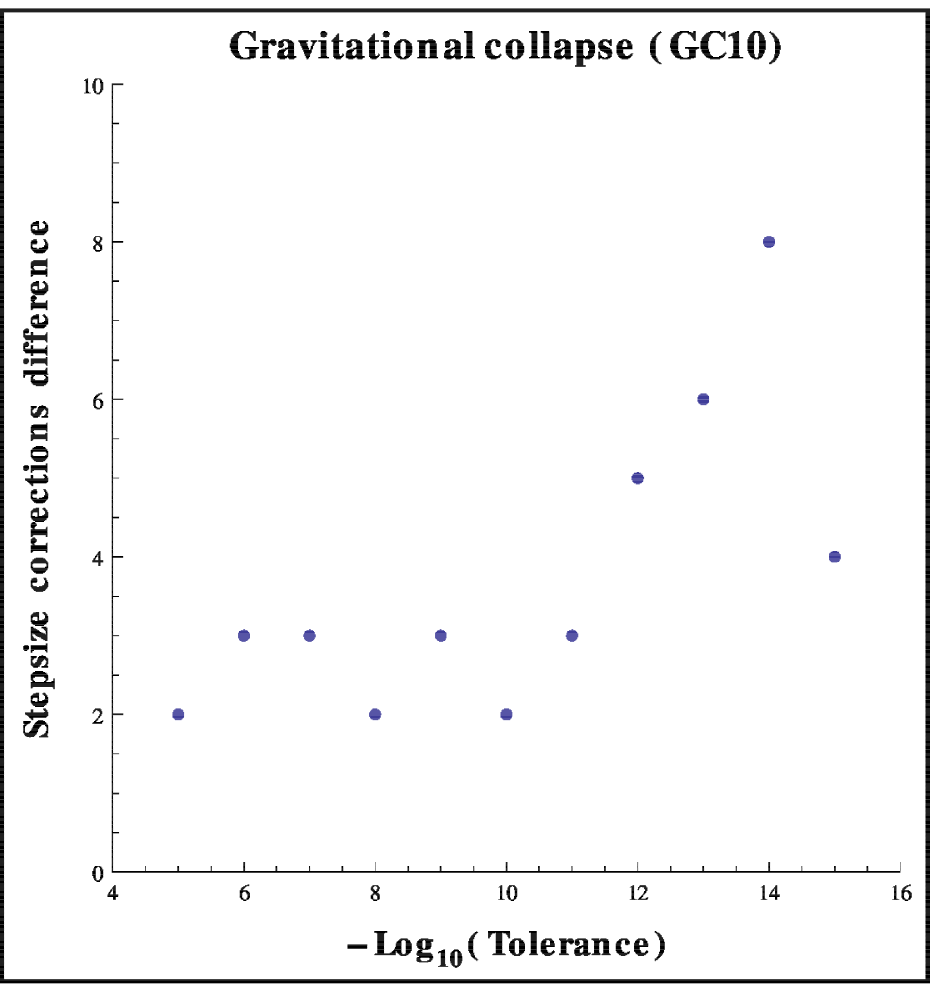}}\quad
\subfigure{\includegraphics[width=3in, height=3in]{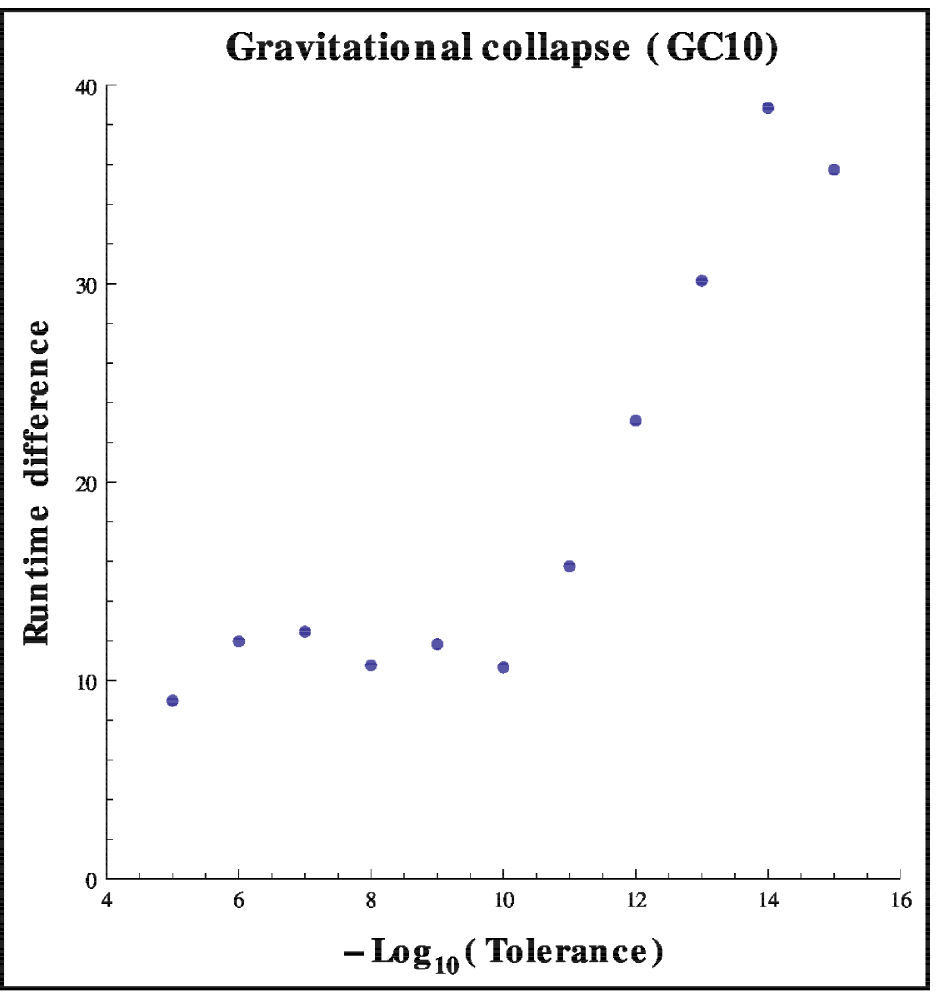}}
}
\mbox{\subfigure{\includegraphics[width=3in, height=3in]{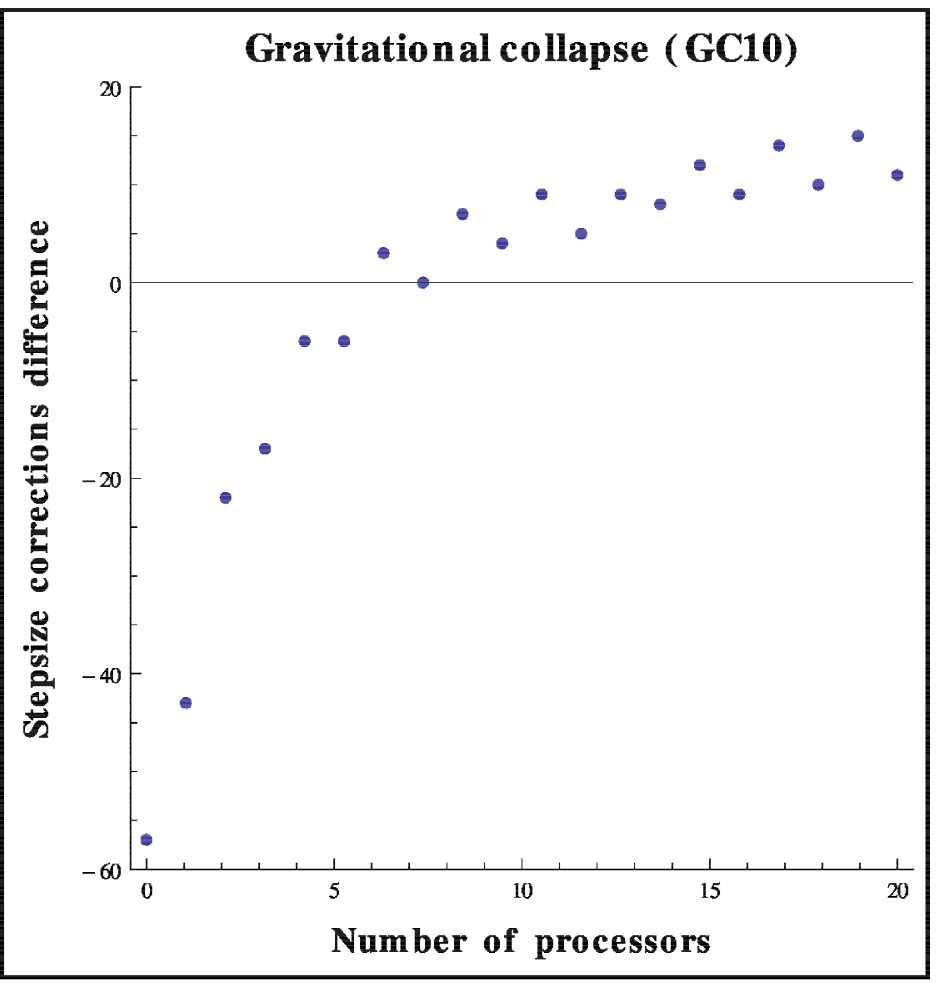}}\quad
\subfigure{\includegraphics[width=3in, height=3in]{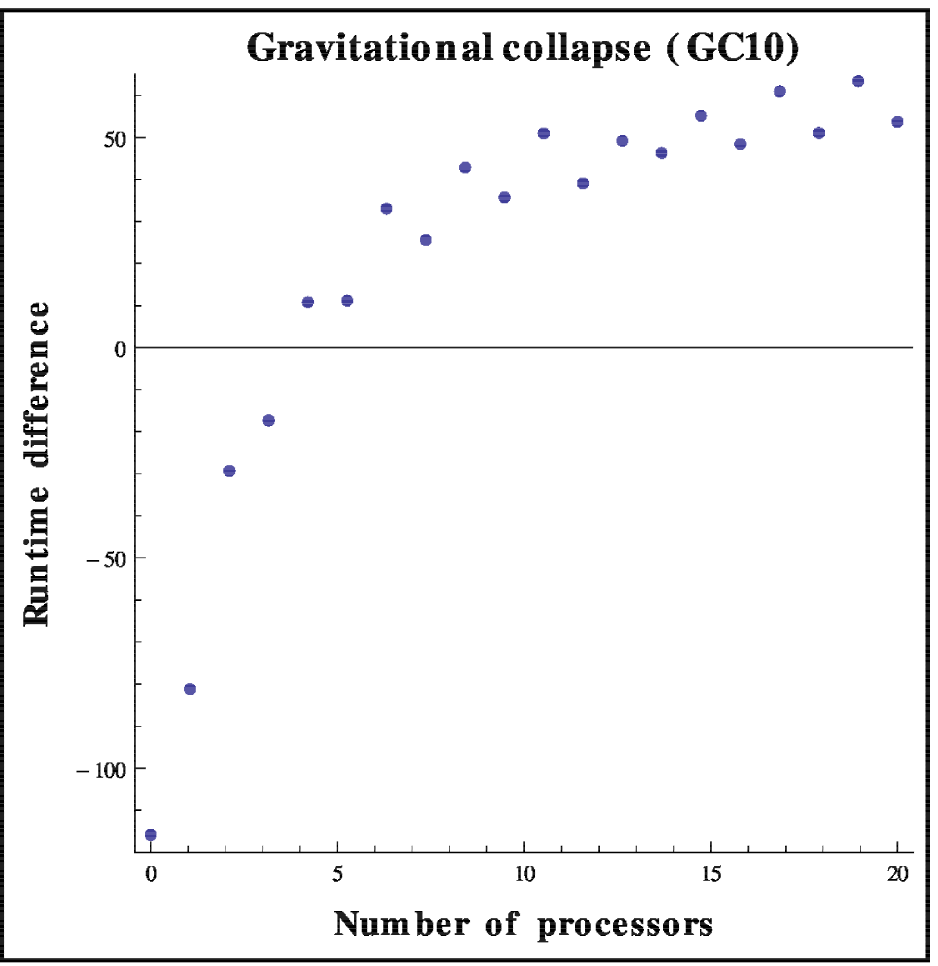}}
}
\caption{DOP853 minus DOP853-ASPA for GC10. Left top: stepsize corrections  
vs tolerance. Left bottom: stepsize corrections vs number of processors. Right 
top: runtime vs tolerance. Right bottom: runtime vs number of processors.} 
\label{fig:ASPAGC10}
\end{figure}

In figures \ref{fig:ASPAGC10} we present the results of the comparison.
For these calculations we used $t_{end}=20$ and $N_{CPU}=20$.
Now it is clear that the parallel algorithm typically performs better than the serial.
From the top two panels we observe that increasing the tolerance induces
a  steadily increase of the difference in required stepsize corrections and 
this 
straightforwardly leads to a larger difference in runtime. The results in the 
bottom panels show that the DOP853-ASPA is efficient when $N_{CPU}>5$. 
Recalling 
that in the bottom panel the differences are all calculated with respect of the 
fixed number obtained with DOP853, we see now that increasing $N_{CPU}$ leads 
to 
less stepsize corrections required by the DOP853-ASPA, but now this also 
corresponds to less runtime. All the above suggests that, indeed, for the GC10 
system, the parallel overhead problem is solved.

To end the comparisons,
note in table \ref{table:cpuvstime} that for this last system, DOP853-ASPA with 
5 processors lasted a little bit more than 7 minutes, while we verified that 
PIRK10 took about an hour to solve the same problem. On the other hand, we also 
checked that for GC40, DOP853-ASPA took about an hour to reach $t_{end}=0.1$ 
with a tolerance of $10^{-6}$, while as mentioned in section \ref{sec:tresults}
PIRK10 needed about six times more (see table \ref{table:colapso1}).

\section{Conclusions}\label{sec:concl}

We tested a parallel iterated Runge--Kutta method of order 10 (PIRK10)
and an adaptive stepsize parallel algorithm (ASPA), 
introduced in this paper, 
which was coupled to a Dormand--Prince method of order 8 (DOP853-ASPA).
The results presented in this paper show that when the initial value problem 
to solve has a simple to evaluate right--hand side (as is the case in more common 
dynamical systems), even in the best case scenarios for the parallel methods,
their performances were only comparable to the corresponding performance of
a serial Dormand-Prince method of order 8 (DOP853). Therefore, taking into 
account code and algorithmic efficiencies, parallel integration seems to not be 
a good practice.

This negative result seems to be due to a parallel overhead problem, i.e.,
the amount of time required to coordinate parallel tasks is larger than the 
time required for evaluating the system right--hand side. We verified that for 
very complex initial value problems or low tolerances the parallel methods can 
outperform DOP853. For instance, such systems arise while using Galerkin 
projection to solve systems of partial differential equations
or when simulating multi--agent systems. 
In these cases, 
it seems to be more efficient to parallelize the search for an optimal stepsize 
for integration than to parallelize the integration scheme. Indeed, our method, 
DOP853-ASPA, consistently outperformed PIRK10 by almost an order of runtime. 
Moreover, even in some cases where DOP853 did a better job than PIRK10, our 
method was able to solve the corresponding initial value problem in less time 
than both these methods.

A nice feature of ASPA is that it does not relies on a given core integrator, 
it can be coupled to any method with a scheme to estimate the local integration 
error. It can even be another parallel method more efficient than the one 
tested 
here.

\section*{Acknowledgments}	
This research was supported by the Sistema Nacional de Investigadores (M\'exico). 
The work of CAT-E was also partially funded by FRABA-UCOL-14-2013 (M\'exico).

\appendix
\section{Butcher tableau}

Butcher tableau for an implicit Runge--Kutta method of order $10$
ref. \cite{Butcher}.
  \begin{table}[!ht]
    \begin{center}
     \small
       $$\begin{array}{c|c c c c c}
         \frac{1}{2}-\omega_{2}&\omega_{1}&\omega_{1}^{\prime}-\omega_{3}+ 
\omega_{4}^{\prime}&\frac{32}{225}-\omega_{5}& 
\omega_{1}^{\prime}-\omega_{3}-\omega_{4}^{\prime}& \omega_{1}-\omega_{6}\\[.25 
cm]
         \frac{1}{2}-\omega_{2}^{\prime}&\omega_{1}-\omega_{3}^{\prime}+ 
\omega_{4}&\omega_{1}^{\prime}&\frac{32}{225}-\omega_{5}^{\prime}&\omega_{1}^{
\prime}-\omega_{6}^{\prime}& \omega_{1}-\omega_{3}^{\prime}-\omega_{4}\\[.25 
cm] 
         \frac{1}{2}&\omega_{1}+\omega_{7}&\omega_{1}^{\prime}+ 
\omega_{7}^{\prime}&\frac{32}{225}&\omega_{1}^{\prime}-\omega_{7}^{\prime}& 
\omega_{1}-\omega_{7}\\[.25 cm] 
         \frac{1}{2}+\omega_{2}^{\prime}&\omega_{1}+\omega_{3}^{\prime}+
         \omega_{4}&\omega_{1}^{\prime}+\omega_{6}^{\prime}& 
\frac{32}{225}+\omega_{5}^{\prime}& \omega_{1}^{\prime}& 
\omega_{1}+\omega_{^3}^{\prime}-\omega_{4}\\[.25 cm]
         \frac{1}{2}+\omega_{2}&\omega_{1}+\omega_{6}&\omega'_{1}+ 
\omega_{3}+\omega_{4}'&\frac{32}{225}+\omega_{5}&\omega'_{1}+\omega_{3}-\omega'_
{4}& \omega_{1}\\[.25 cm]
         \hline\\[-.25 cm]
          & 2\omega_{1}&2\omega_{1}'&\frac{64}{225}&2\omega_{1}'& 2\omega_{1}\\
          \end{array}$$
            \caption{5-stage, order 10.} \label{table:RK10}
        \end{center}
    \end{table}
    
where the $\omega_i$ are given by,
\begin{align}
&\omega_{1}=\frac{322-13\sqrt{70}}{3600},\quad
\omega_{1}'=\frac{322+13\sqrt{70}}{3600},\nonumber\\
&\omega{2}=\frac{1}{2}\sqrt{\frac{35+2\sqrt{70}}{63}},\quad 
\omega{2}'=\frac{1}{2}\sqrt{\frac{35-2\sqrt{70}}{63}},\nonumber\\
&\omega{3}=\omega_2\frac{452+59\sqrt{70}}{3240},\quad
\omega'_3=\omega'_2\frac{452-59\sqrt{70}}{3240},\nonumber\\
&\omega_4=\omega_2\frac{64+11\sqrt{70}}{1080},\quad
\omega'_4=\omega'_2\frac{64-11\sqrt{70}}{1080},\nonumber\\
&\omega_5=8\omega_2\frac{23-\sqrt{70}}{405},\quad
\omega'_5=8\omega'_2\frac{23+\sqrt{70}}{405},\nonumber\\
&\omega_6=\omega_2-2\omega_3-\omega_5,\quad
\omega'_6=\omega'_2-2\omega'_3-\omega'_5,\nonumber\\
&\omega_7=\omega_2\frac{308-23\sqrt{70}}{960},\quad 
\omega'_7=\omega'_2\frac{308+23\sqrt{70}}{960}.\nonumber
\end{align}

\section{Runtimes}
\label{sec:runtimes}
  \begin{table}[ht]
    \begin{center}
    \begin{tabular}{c|c|c|c|c}
     & \multicolumn{2}{|c|}{DOP853-ASPA} &\multicolumn{2}{|c}{DOP853}\\
     \hline
    T   & stepsize corrections& Time & stepsize corrections & Time\\
    \hline
    5   &  4333 & 0.06 &  4641 & 0.09 \\[.2 cm]
    6   &  5778 & 0.07 &  6155 & 0.02 \\[.2 cm]
    7   &  7826 & 0.1  &  8186 & 0.03 \\[.2 cm]
    8   & 10520 & 0.12 & 10846 & 0.03 \\[.2 cm]
    9   & 14040 & 0.17 & 14394 & 0.04 \\[.2 cm]
    10  & 18825 & 0.27 & 18928 & 0.05 \\[.2 cm]
    11  & 24876 & 0.29 & 25297 & 0.07 \\[.2 cm]
    12  & 33327 & 0.38 & 33657 & 0.09 \\[.2 cm]
    13  & 49989 & 0.52 & 45424 & 0.12 \\[.2 cm]
    14  & 59292 & 0.67 & 60145 & 0.12 \\[.2 cm]
    15  & 79393 & 0.91 & 79990 & 0.21 \\
    \hline
    \end{tabular}
    \caption{$T=-\log_{10}(tolerance)$. HH system using 10 processor and final 
time 5000.}
    \end{center}
  \end{table}

      \begin{table}[ht]
    \begin{center}
    \begin{tabular}{c|c|c|c|c}
     & \multicolumn{2}{|c|}{DOP853-ASPA} &\multicolumn{2}{|c}{DOP853}\\
    \hline
    CPU's & stepsize corrections & Time & stepsize corrections & Time \\
    \hline
    1   & 131089& 0.53 & 79990 & 0.21 \\[.2 cm]
    2   & 97312 & 0.59 & - & -\\[.2 cm]
    3   & 89886 & 0.6  & - & -\\[.2 cm]
    4   & 87133 & 0.67 & - & -\\[.2 cm]
    5   & 82941 & 0.7  & - & -\\[.2 cm]
    6   & 83143 & 0.72 & - & -\\[.2 cm]
    7   & 79923 & 0.78 & - & -\\[.2 cm]
    8   & 80891 & 0.84 & - & -\\[.2 cm]
    9   & 78640 & 0.83 & - & -\\[.2 cm]
    10  & 79393 & 0.93 & - & -\\[.2 cm]
    11  & 77713 & 0.96 & - & -\\[.2 cm]
    12  & 77428 & 0.99 & - & -\\[.2 cm]
    13  & 77469 & 1.04 & - & -\\[.2 cm]
    14  & 77170 & 1.05 & - & -\\[.2 cm]
    15  & 75664 & 1.1  & - & -\\[.2 cm]
    16  & 76341 & 1.19 & - & -\\[.2 cm]
    17  & 76243 & 1.33 & - & -\\[.2 cm]
    18  & 76566 & 1.43 & - & -\\[.2 cm]
    19  & 75997 & 1.54 & - & -\\[.2 cm]
    20  & 75799 & 1.8  & - & -\\
    \hline
    \end{tabular}
    \caption{HH system using tolerance $10^{-15}$.}
    \end{center}
  \end{table}

  \begin{table}[ht]
    \begin{center}
    \begin{tabular}{c|c|c|c|c}
     & \multicolumn{2}{|c|}{DOP853-ASPA} &\multicolumn{2}{|c}{DOP853}\\
     \hline
    T   & stepsize corrections& Time & stepsize corrections & Time\\
    5   & 4321 & 0.74 & 4662  & 0.75\\[.2 cm]
    6   & 5853 & 1.01 & 6204  & 0.99\\[.2 cm]
    7   & 7844 & 1.34 & 8170  & 1.26\\[.2 cm]
    8   & 10499& 1.79 & 10814 & 1.62\\[.2 cm]
    9   & 13829& 2.35 & 14291 & 2.03\\[.2 cm]
    10  & 18660& 3.16 & 19086 & 2.61\\[.2 cm]
    11  & 24902& 4.22 & 25422 & 3.57\\[.2 cm]
    12  & 33369& 5.66 & 33740 & 4.53\\[.2 cm]
    13  & 44362& 7.53 & 45267 & 6.09\\[.2 cm]
    14  & 59127& 10.07& 60171 & 8.27\\[.2 cm]
    15  & 79444& 13.45& 80568 & 10.96\\
    \hline
    \end{tabular}
    \caption{$T=-\log_{10}(tolerance)$. HH100 system. Here we used 10 
processors and final time 5000.}
    \end{center}
  \end{table}

    \begin{table}[H]
    \begin{center}
    \begin{tabular}{c|c|c|c|c}
    & \multicolumn{2}{|c|}{DOP853-ASPA} &\multicolumn{2}{|c}{DOP853}\\
    \hline
    CPU's & stepsize corrections & Time & stepsize corrections & Time \\
    \hline
    1   & 131089& 17.79 & 80568 & 10.96\\[.2 cm]
    2   & 97765 & 14.02 & - & -\\[.2 cm]
    3   & 89631 & 13.16 & - & -\\[.2 cm]
    4   & 87586 & 13.09 & - & -\\[.2 cm]
    5   & 83745 & 12.81 & - & -\\[.2 cm]
    6   & 83174 & 12.97 & - & -\\[.2 cm]
    7   & 80635 & 12.81 & - & -\\[.2 cm]
    8   & 80946 & 13.18 & - & -\\[.2 cm]
    9   & 78234 & 13.08 & - & -\\[.2 cm]
    10  & 79444 & 13.48 & - & -\\[.2 cm]
    11  & 77802 & 13.33 & - & -\\[.2 cm]
    12  & 78105 & 13.88 & - & -\\[.2 cm]
    13  & 76761 & 13.67 & - & -\\[.2 cm]
    14  & 77358 & 14.04 & - & -\\[.2 cm]
    15  & 77139 & 14.25 & - & -\\[.2 cm]
    16  & 76813 & 14.56 & - & -\\[.2 cm]
    17  & 76853 & 14.86 & - & -\\[.2 cm]
    18  & 76150 & 14.97 & - & -\\[.2 cm]
    19  & 76295 & 17.33 & - & -\\[.2 cm]
    20  & 76440 & 15.59 & - & -\\
    \hline
    \end{tabular}
     \caption{HH100 system.}
    \end{center}
  \end{table}
  
  \begin{table}[H]
    \begin{center}
    \begin{tabular}{c|c|c|c|c}
     & \multicolumn{2}{|c|}{DOP853-ASPA} &\multicolumn{2}{|c}{DOP853}\\
     \hline
    T   & stepsize corrections& Time & stepsize corrections & Time\\
    \hline
    5   & 7  & 17.39 & 9  & 26.37\\[.2 cm]
    6   & 10 & 24.82 & 13 & 36.79\\[.2 cm]
    7   & 15 & 37.25 & 18 & 49.71\\[.2 cm]
    8   & 22 & 54.65 & 24 & 65.42\\[.2 cm]
    9   & 30 & 74.45 & 33 & 86.28\\[.2 cm]
    10  & 42 & 104.37& 44 & 115.03\\[.2 cm]
    11  & 57 & 141.7 & 60 & 157.45\\[.2 cm]
    12  & 75 & 186.07& 80 & 209.17\\[.2 cm]
    13  & 102& 253.1 & 108& 283.24\\[.2 cm]
    14  & 136& 337.85& 144& 376.71\\[.2 cm]
    15  & 189& 469.16& 193& 504.9 \\
    \hline
    \end{tabular}
    \caption{$T=-\log_{10}(tolerance)$. Collapse (GC10) using 10 processors 
    and final time 20.}
    \end{center}
\end{table} 

  \begin{table}[H]
    \begin{center}
    \begin{tabular}{c|c|c|c|c}
    & \multicolumn{2}{|c|}{DOP853-ASPA} &\multicolumn{2}{|c}{DOP853}\\
    \hline
    CPU's & stepsize corrections & Time & stepsize corrections & Time \\
    \hline
    1   & 250 & 620.75 & 193 & 504.9\\[.2 cm]
    2   & 236 & 586.12 & - & -\\[.2 cm]
    3   & 215 & 534.19 & - & -\\[.2 cm]
    4   & 210 & 522.19 & - & -\\[.2 cm]
    5   & 199 & 494.1  & - & -\\[.2 cm]
    6   & 199 & 493.71 & - & -\\[.2 cm]
    7   & 190 & 471.84 & - & -\\[.2 cm]
    8   & 193 & 479.29 & - & -\\[.2 cm]
    9   & 186 & 462.08 & - & -\\[.2 cm]
    10  & 189 & 469.16 & - & -\\[.2 cm]
    11  & 184 & 453.94 & - & -\\[.2 cm]
    12  & 188 & 465.82 & - & -\\[.2 cm]
    13  & 184 & 455.71 & - & -\\[.2 cm]
    14  & 185 & 458.62 & - & -\\[.2 cm]
    15  & 181 & 449.76 & - & -\\[.2 cm]
    16  & 184 & 456.48 & - & -\\[.2 cm]
    17  & 179 & 443.96 & - & -\\[.2 cm]
    18  & 183 & 453.85 & - & -\\[.2 cm]
    19  & 178 & 441.52 & - & -\\[.2 cm]
    20  & 182 & 451.2  & - & -\\
    \hline
    \end{tabular}
    \caption{GC10 using tolerance $10^{-15}$ and final time 20.} 
\label{table:cpuvstime}
    \end{center}
  \end{table}

\end{document}